\documentclass[10pt,twoside]{classe-ma_e}
\usepackage{amssymb,amsbsy,amsmath,amsfonts,amssymb,amscd}
\usepackage{latexsym,euscript,exscale,epsfig}
\usepackage[english,francais]{babel}
\usepackage{times}
\newtheorem{e-proposition}[theorem]{Proposition}

\newtheorem{e-definition}[theorem]{Definition\rm}
\newtheorem{remark}{\it Remark\/}

\newtheorem{notation}{\it Notation\/}

\ComParit{S0764-4442}
\PIT{FLA}
\PXMA{????}
\Add{?}
\Volume{???}
\Year{2001}
\FirstPage{1}
\LastPage{6}
\AuteurCourant{Matthieu Picantin}
\TitreCourant{Presentations for the Dual Braid Monoids} 
\Journal
\Rubrique{Th\'eorie des groupes}{Group Theory}
\SousRubrique{}{}
\PresentePar{First name}{NAME}
\Recu{jour mois ann\'ee}{apr\`es r\'evision}{jour mois ann\'ee}
\begin{document}
\selectlanguage{english}
\title{%
Explicit Presentations for the Dual Braid Monoids
}
\author{%
Matthieu Picantin\ \
}
\address{%
Laboratoire SDAD, D\'epartement de Math\'ematiques, Universit\'e de Caen, 14000 Caen, France \\
E-mail: picantin@math.unicaen.fr
}
\maketitle
\thispagestyle{empty}
\begin{Abstract}{%
Birman, Ko \&~Lee have introduced a new monoid~${\cal B}^{\!*\!}_{\!n}$---with an explicit presentation---whose group
of fractions is the \mbox{$n$-strand} braid group~${\cal B}_{\!n}$. Building on a new approach by~Digne, Michel and~himself, Bessis
has defined a {\it dual} braid monoid for every finite Coxeter type Artin-Tits group extending the type~A case. Here, we give
an explicit presentation for this dual braid monoid in the case of types~B and~D, and we study the combinatorics of the
underlying Garside structures.}
\end{Abstract}
\selectlanguage{french}
\begin{Ftitle}{%
Pr\'esentations pour les mono\"\i des de tresses duaux
}\end{Ftitle}
\begin{Resume}{%
Birman, Ko et~Lee ont introduit un nouveau mono\"\i de~${\cal B}^{\!*\!}_{\!n}$---avec une pr\'esentation
explicite---dont le groupe de fractions est le groupe~${\cal B}_{\!n}$ des tresses \`a~$n$ brins. Suivant une nouvelle
approche propos\'ee avec~Digne et~Michel, Bessis a d\'efini un mono\"\i de de tresses {\rm dual} pour tout groupe d'Artin-Tits
de type de Coxeter fini g\'en\'eralisant le cas du type~A. Ici, nous donnons une pr\'esentation explicite de ce mono\"\i de de
tresses dual pour les groupes d'Artin-Tits de type~B et~D, et nous \'etudions la combinatoire des structures de Garside
sous-jacentes.}
\end{Resume}
\AFv

\font\ttX=cmtt10 at 9pt
\def\ie{\hbox{\it i.e.}}
\def\ni{\noindent}
\def\bgni{\bigbreak\goodbreak\noindent}
\def\bnni{\bigbreak\nobreak\noindent}
\def\mgni{\medbreak\goodbreak\noindent}
\def\mnni{\medbreak\nobreak\noindent}
\def\sgni{\smallbreak\goodbreak\noindent}
\def\snni{\smallbreak\nobreak\noindent}

\def\scr{\scriptscriptstyle}
\def\BB{{\mathbf{B}}}
\def\gart{\hbox{$\BB$}}
\def\clas{\hbox{$\BB^{\scr\!+\!}$}}
\def\dual{\hbox{$\BB^{\!*\!}$}}
\def\vsm{\vspace{-0,2cm}}
\def\a{\alpha}
\def\b{\beta}
\def\d{\delta}
\def\D{\Delta}
\def\g{\gamma}
\def\s{\sigma}
\def\t{\tau}
\ni Birman, Ko \&~Lee introduisent dans~\cite{bkl} un nouveau mono\"\i de pour
les groupes de tresses (de type~A) avec une pr\'esentation explicite. La question de possibles g\'en\'eralisations se pose
naturellement. Une bonne notion pour l'\'etude de ces nouveaux mono\"\i des de tresses (mais aussi de mono\"\i des pour les
groupes de tresses des groupes de r\'eflexions complexes, pour les groupes d'entrelacs, etc) est celle de mono\"\i de
de~Garside, introduite par Dehornoy \&~Paris dans~\cite{dfx} et exploit\'ee dans~\cite{bes,bdm,dhh,dgk,pia,pib,pic,pie}~: $M$
est un {\it mono\"\i de de~Garside} si~$M$ est un mono\"\i de simplifiable, admet des ppcm \`a droite et \`a gauche et admet
un {\it \'el\'ement de~Garside} d\'efini comme un \'el\'ement dont les diviseurs \`a droite et \`a gauche co\"\i
ncident, engendrent~$M$ et sont en nombre fini. Les diviseurs de l'\'el\'ement de~Garside minimal sont
appel\'es les {\it
\'el\'ements simples} de~$M$~; muni des op\'erations ppcm et pgcd, l'ensemble des \'el\'ements simples est un treillis fini.
Les mono\"\i des de~Garside se plongent dans leurs groupes de fractions, ont de bonnes formes normales,
des structures automatiques explicites, etc. Le crit\`ere donn\'e dans~\cite{dgk} permet de d\'ecider si une pr\'esentation de
mono\"\i de est celle d'un mono\"\i de de~Garside~: il consiste en la v\'erification de conditions de {\it
compl\'etude} et de~{\it cube} et de l'existence d'un \'el\'ement de~Garside.

Dans une pr\'esentation (de mono\"\i de), si~$w_{\scr\!1},\ldots,w_{\scr\!p}$ sont des mots, nous
\'ecrivons~$[w_{\scr\!1},\ldots,w_{\scr\!p}]$ pour la famille de
relations~$w_{\scr\!1}w_{\scr\!2}=w_{\scr\!2}w_{\scr\!3}=\ldots=w_{\scr\!p-1}w_{\scr\!p}=w_{\scr\!p}w_{\scr\!1}$
(compatible avec le symbole du commutateur dans une pr\'esentation de groupe).

Le principal r\'esultat de~\cite{bkl} peut s'\'enoncer comme suit~: le sous-mono\"\i de~$\dual(A_{n-1})$ du groupe de tresses
d'Artin-Tits~$\gart(A_{n-1})$ engendr\'e par les
tresses~$a_{ts}=(\s_{t-1}\cdots\s_{s+1})\s_s(\s_{t-1}\cdots\s_{s+1})^{\scr-1}$ pour~$n\geq t>s\geq1$ (o\`u les
$\s_i$ sont les g\'en\'erateurs du mono\"\i de classique~$\clas(A_{n-1})$) admet la pr\'esentation~(\ref{eq:presa})~: c'est un
mono\"\i de de~Garside, dont le nombre de simples est le $n$-i\`eme nombre de~Catalan (voir la table~\ref{ta:table}).

Suivant une nouvelle approche propos\'ee avec~Digne \&~Michel dans~\cite{bdm}, en g\'en\'eralisant le cas
du type~A, Bessis d\'efinit dans~\cite{bes} un mono\"\i de de tresses {\it dual}~$\dual(T)$
pour tout groupe d'Artin-Tits~$\gart(T)$ de type de Coxeter fini~T, comme \'etant {\it le} mono\"\i de
de~Garside dont le treillis des simples est isomorphe au treillis de~$\prec$-divisibilit\'e d'un
\'el\'ement de Coxeter dans le groupe de Coxeter associ\'e, o\`u~$\prec$ est d\'efinie relativement \`a la
longueur en r\'eflexions (les preuves pour le cas du type~D ne sont pas encore publi\'ees).

Nous montrons dans cette note que le mono\"\i de de tresses dual~$\dual(B_n)$ admet la pr\'esentation~(\ref{eq:presb}), tandis que le
mono\"\i de de tresses dual~$\dual(D_n)$ admet la pr\'esentation~(\ref{eq:presd}). Les preuves consistent \`a montrer que le
sous-mono\"\i de du groupe~$\gart(B_n)$ (resp. $\gart(D_n)$) engendr\'e par les g\'en\'erateurs d\'efinis par~(\ref{eq:genb})
(resp. par~(\ref{eq:gend})) admet la pr\'esentation~(\ref{eq:presb}) (resp. (\ref{eq:presd})) en utilisant les diagrammes de
tresses des figures~\ref{fi:genb} et~\ref{fi:bisotopy} (resp. \ref{fi:gend} et~\ref{fi:disotopy}), et que cette
pr\'esentation est celle d'un mono\"\i de de~Garside, dont l'\'el\'ement de Garside minimal a pour image un \'el\'ement
de~Coxeter dans le groupe de~Coxeter associ\'e.

Une approche analogue \`a celle de Birman, Ko \&~Lee~\cite{bkl} permet de montrer que les
\'el\'ements simples de~$\dual(B_n)$ et~$\dual(D_n)$ sont en bijection avec les partitions non-crois\'ees
correspondantes que~Reiner d\'efinit dans~\cite{rei}. La table~\ref{ta:table} donne le nombre
d'\'el\'ements simples pour les mono\"\i des de tresses duaux, rassemblant les r\'esultats th\'eoriques
pour les types~A, B, D, I$_2$ et des r\'esultats obtenus par le calcul---utilisant le progiciel~CHEVIE
de~GAP~\cite{gap}---pour les types exceptionnels (le r\'esultat du calcul pour~E$_8$ appara\^\i t d\'ej\`a
dans~\cite{bes}).

\par\medskip\centerline{\rule{2cm}{0.2mm}}\medskip
\setcounter{section}{0}
\selectlanguage{english}

\section{Introduction}\label{intro}

\ni Birman, Ko \&~Lee introduced in~\cite{bkl} an alternative monoid for braid groups (of
type~A) together with an explicit presentation of this monoid. The question of possible generalizations arises naturally. A good
notion for studying such new braid monoids (but also monoids for braid groups of complex reflection groups, for link groups,
etc) is that of a Garside monoid, introduced by Dehornoy \&~Paris in~\cite{dfx} and further studied
in~\cite{bes,bdm,dhh,dgk,pia,pib,pic,pie}~: $M$ is a {\it Garside monoid} if~$M$ is cancellative,
admits right and left lcm's and admits a {\it Garside element} defined to be an element whose left and right
divisors coincide, generate~$M$ and are finite in number. The divisors of the minimal Garside element are
called {\it simple elements} of~$M$~; when equipped with lcm and gcd operations, the set of simple
elements is a finite lattice. Garside monoids embed into their groups of fractions, they admit nice normal
forms, explicit automatic structures, etc. Whether a given monoid presentation is that a Garside monoid
can be decided using Dehornoy's criterion of~\cite{dgk}~: it consists in the verification of some {\it
completeness} and~{\it cube} conditions and of the existence of a Garside element.

\begin{notation}In a (monoid) presentation, $w_{\scr\!1},\ldots,w_{\scr\!p}$ being words, we write~$[w_{\scr\!1},\ldots,w_{\scr\!p}]$
for~$w_{\scr\!1}w_{\scr\!2}=w_{\scr\!2}w_{\scr\!3}=\ldots=w_{\scr\!p-1}w_{\scr\!p}=w_{\scr\!p}w_{\scr\!1}$ (which is compatible with the {\it
commutator} notation in a group presentation).
\end{notation}

\begin{e-proposition}{\rm\cite{bkl}} \label{presa}The submonoid~$\dual(A_{n-1})$ of the Artin-Tits braid group~$\gart(A_{n-1})$
generated by~$a_{ts}=(\s_{t-1}\cdots\s_{s+1})\s_s(\s_{t-1}\cdots\s_{s+1})^{\scr-1}$ for~$n\geq t>s\geq1$ (where the
$\s_i$'s are the generators for the classical monoid~$\clas(A_{n-1})$) admits the presentation
\begin{equation}\label{eq:presa}
  \langle~a_{ts}:~[a_{ts},a_{sr},a_{tr}]~~\hbox{for}~t>s>r~,
  ~[a_{ts},a_{rq}]~~\hbox{for}~(t\!-\!r)(t\!-\!q)(s\!-\!r)(s\!-\!q)>0~\rangle;
\end{equation}it is a Garside monoid, whose number of simple elements is the $n$-th
Catalan number (see Table~\ref{ta:table}).
\end{e-proposition}

\begin{e-proposition}\label{presi}The submonoid of~$\gart(I_2(m))$ generated
by~$\s_1$ and~$\s_i=(\overbrace{\s_2\s_1\s_2\cdots}^{\scr(i-1)\rm{~terms}})
{(\overbrace{\s_2\s_1\s_2\cdots}^{\scr(i-2)\rm{~terms}}})^{\scr-1}$ for~$2\leq i\leq m$ (where the
$\s_i$'s are the generators for the classical monoid~$\clas(I_2(m))$) is presented
by~$\langle~\s_i~:~[\s_m,\ldots,\s_1]~\rangle$;  it is a Garside monoid, denoted by~$\dual(I_2(m))$.
\end{e-proposition}

\mnni Building on a new approach by~Digne, Michel and himself in~\cite{bdm},
Bessis defined in~\cite{bes}, extending the type~A case, a {\it dual} braid monoid~$\dual(T)$ for every
finite Coxeter type~T Artin-Tits group~$\gart(T)$ to be the Garside monoid whose lattice of simple
elements is the
$\prec$-divisibility lattice of a Coxeter element in the associated Coxeter group, where~$\prec$ is
defined with respect to the reflection length (proofs for the type~D case are not published yet).
\section{Dual monoids for type~B Artin-Tits braid groups}\label{btype}

\begin{figure}[!htb]
\begin{center}
\scalebox{.90}{\includegraphics[width=12cm,height=2cm]{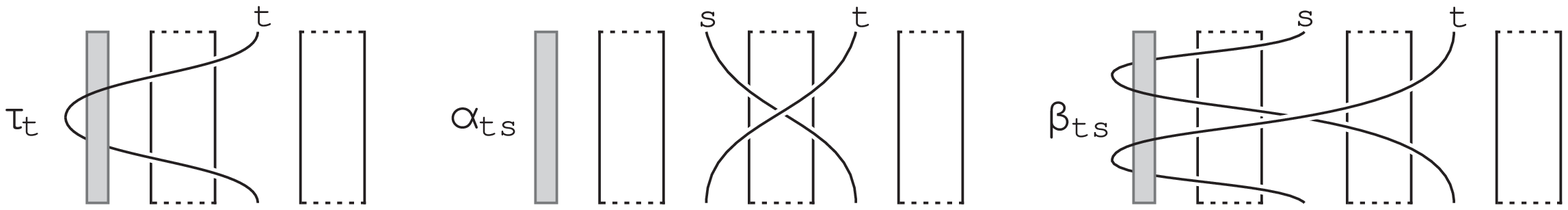}}
\end{center}
\vglue -0.3cm
\caption{Braid pictures for type~B new generators.}\label{fi:genb}
\centerline{\it Figure~\ref{fi:genb}: Diagrammes de tresse pour les nouveaux g\'en\'erateurs de type~B.}
\end{figure}

\ni In this section, we establish analogous to Propositions~\ref{presa} and~\ref{presi} for type~B. The classical monoid~$\clas(B_n)$ for the
Artin-Tits braid group~$\gart(B_n)$ admits the presentation\begin{equation}\label{eq:presclasb}
 \begin{aligned}
  \langle~\t_1^{},\s_1^{},\ldots,\s_{n-1}^{}:
   &~\s_1^{}\t_1^{}\s_1^{}\t_1^{}=\t_1^{}\s_1^{}\t_1^{}\s_1^{}~,
   ~\s_i^{}\s_{i+1}^{}\s_i^{}=\s_{i+1}^{}\s_i^{}\s_{i+1}^{}~,~1\leq i\leq n-2\\
   &~\s_i^{}\s_j^{}=\s_j^{}\s_i^{}~,~1<i+1<j<n~,
   ~\t_1^{}\s_j^{}=\s_j^{}\t_1^{}~,~1<j<n~\rangle.
 \end{aligned}
\end{equation}We shall use the well-known fact that~$\gart(B_n)$
can be viewed as the subgroup of~$\gart(A_n)$ of those braids
whose first strand is not braided. Let us introduce the
following~$n^2$ new generators~:\begin{equation}\label{eq:genb}
 \begin{aligned}\a_{ts}^{}&=(\s_{t-1}^{}\s_{t-2}^{}\cdots\s_{s+1}^{})
  \s_s^{}(\s_{t-1}^{}\s_{t-2}^{}\cdots\s_{s+1}^{})^{-1}\quad\hbox{for}\quad n\geq t>s\geq 1,\\
  \t_1 \hbox{ and }\t_t^{}&=\a_{t1}^{}\t_1^{}\a_{t1}^{-1}\quad\hbox{for}\quad n\geq t>1,\\
  \b_{ts}^{}&=\t_s^{-1}\a_{ts}^{}\t_s^{}\quad\hbox{for}\quad n\geq t>s\geq 1.
 \end{aligned}
\end{equation}Braid pictures for new generators are displayed in
Figure~\ref{fi:genb} (a ribbon indicates some number of strands
moving in parallel, making some pictures easier to understand).

\begin{figure}[!htb]
\begin{center}
\scalebox{.90}{\includegraphics[width=15cm,height=4cm]{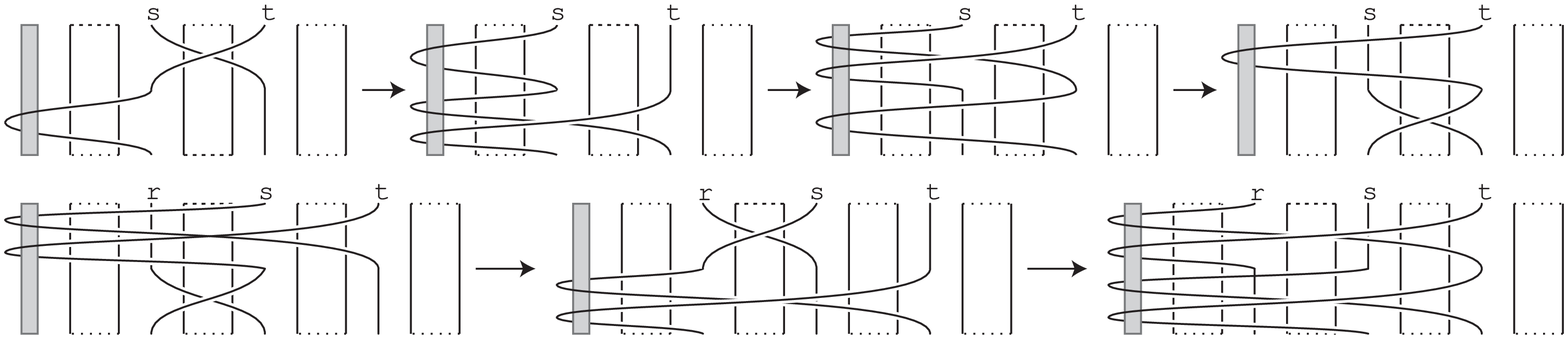}}
\end{center}
\vglue -0.3cm
\caption{Type~B braid isotopies.}\label{fi:bisotopy}
\centerline{\it Figure~\ref{fi:bisotopy}: Isotopies de tresses de type~B.}
\end{figure}

\begin{e-proposition}\label{presb} The dual braid monoid~$\dual(B_n)$ admits the
presentation\begin{equation}\label{eq:presb}
 \begin{aligned}\langle~\a_{ts},\b_{ts},\t_t:
  &~[\a_{ts},\t_s,\b_{ts},\t_t]~~\hbox{for}~t>s,\\
  &~[\a_{ts},\a_{sr},\a_{tr}]~,
   ~[\b_{ts},\a_{sr},\b_{tr}]~,
   ~[\a_{ts},\b_{sr},\b_{tr}]~~\hbox{for}~t>s>r,\\
  &~[\a_{ts},\t_{r}]~,
   ~[\t_t,\a_{sr}]~,
   ~[\b_{tr},\t_{s}]~~\hbox{for}~t>s>r,\\
  &~[\a_{ts},\a_{rq}]~,
   ~[\a_{ts},\b_{rq}]~,
   ~[\b_{ts},\a_{rq}]~,
   ~[\a_{tq},\a_{sr}]~,
   ~[\b_{tq},\a_{sr}]~,
   ~[\b_{tq},\b_{sr}]~~\hbox{for}~t>s>r>q~\rangle.
 \raisetag{1.5cm}\end{aligned}
\end{equation}
\end{e-proposition}
\begin{proof}We first show that the submonoid of~$\gart(B_n)$
generated by the generators of~(\ref{eq:genb}) admits
Presentation~(\ref{eq:presb}), and then that this submonoid is a
Garside monoid by using the criterion given in~\cite{dgk}. 

(i) The type~B braid isotopies displayed in Figure~\ref{fi:bisotopy} give on
the one hand~$[\a_{ts},\t_s,\b_{ts},\t_t]$ for~$n\geq t>s\geq 1$ and on the
other hand~$[\b_{ts},\a_{sr},\b_{tr}]$ for~$t>s>r$. All relations in~(\ref{eq:presb}) are obtained
similarly. Conversely, not so tedious computations prove that the relations of~(\ref{eq:presclasb}) are
consequences of those of~(\ref{eq:presb}).

(ii) The completion algorithm of~\cite{djj} can be successfully applied to~Presentation~(\ref{eq:presb}).
We obtain~:\begin{equation}
\begin{align*}
\b_{sr}^{}\t_s^{}\b_{tr}^{}&=\a_{tr}^{}\a_{ts}^{}\t_r^{}
=\b_{ts}^{}\a_{sr}^{}\t_t^{}=\t_t^{}\a_{sr}^{}\a_{tr}^{}=\t_s^{}\b_{ts}^{}\a_{sr}^{}
=\b_{tr}^{}\t_t^{}\a_{ts}^{}=\t_s^{}\b_{ts}^{}\a_{sr}^{},
\end{align*}
\end{equation}
for~$t>s>r$, and\begin{equation*}
\begin{array}{rclrclrcl}
 \a_{tq}^{}\a_{ts}^{}\a_{sr}^{}\t_q^{}\!\!&\!\!=\!\!&\!\!\b_{sr}^{}\b_{tr}^{}\a_{rq}^{}\t_s^{},
 &\b_{ts}^{}\a_{sq}^{}\a_{sr}^{}\t_t^{}\!\!&\!\!=\!\!&\!\!\b_{rq}^{}\b_{tq}^{}\a_{ts}^{}\t_r^{},
 &\a_{tr}^{}\a_{rq}^{}\a_{ts}^{}\!\!&\!\!=\!\!&\!\!\a_{sq}^{}\a_{sr}^{}\a_{tq}^{},\\
 \a_{tr}^{}\a_{ts}^{}\b_{rq}^{}\!\!&\!\!=\!\!&\!\!\b_{sq}^{}\a_{sr}^{}\b_{tq}^{},
 &\b_{tr}^{}\a_{rq}^{}\b_{ts}^{}\!\!&\!\!=\!\!&\!\!\a_{sq}^{}\a_{sr}^{}\b_{tq}^{},
 &\b_{tr}^{}\a_{rq}^{}\a_{ts}^{}\!\!&\!\!=\!\!&\!\!\b_{sq}^{}\b_{sr}^{}\b_{tq}^{},
\end{array}
\end{equation*}
for~$t>s>r>q$. Now, every relation in the complete presentation involves at most 4 distinct strands, so, in order to prove
local cube condition, hence, by homogeneity, global cube condition, it suffices to check local cube condition of the complete presentation
obtained for~$\dual(B_5)$---which is immediate with a computer. Finally,
we verify that~$\d=\a_{n,n-1}\cdots\a_{2,1}\t_1$ is a Garside element, whose image in the associated
Coxeter group is a Coxeter element.
\end{proof}

\begin{remark} Another proof of~Proposition~\ref{presb}(ii) is as follows. From~\cite[Theorems~9.2, 9.3
and~9.4]{dfx}, we deduce that the monoid~$\dual(A_{2n-1})^\phi$ of elements fixed under the halfturn
automorphism~$\phi$ of~$\dual(A_{2n-1})$ is a Garside monoid. If the atoms of~$\dual(A_{2n-1})$ are denoted
by~$a_{ts}$ for~$2n\geq t>s\geq1$ (see Prop.~\ref{presa}), the atoms of~$\dual(A_{2n-1})^\phi$ are
$a_{(n+i)i}$, $a_{(n+j+i-1)(n+i)}a_{(j+i-1)i}$ for~$n\!\geq\!i\!\geq\!1$ \mbox{and~$n\!\geq\!j\!\geq\!2$} (with
indices taken modulo~$2n$ and~$a_{st}=a_{ts}$). Now, the map defined \mbox{by~$\tau_1\mapsto
a_{(n+1)1}$} and~$\a_{(i+1)i}\mapsto a_{(n+1+i)(n+i)}a_{(i+1)i}$ extends into an isomorphism
from~$\dual(B_n)$ to~$\dual(A_{2n-1})^\phi$.
\end{remark}

\section{Dual monoids for type~D Artin-Tits braid groups}\label{dtype}

\begin{figure}[!htb]
\begin{center}
\scalebox{.90}{\includegraphics[width=16cm,height=2cm]{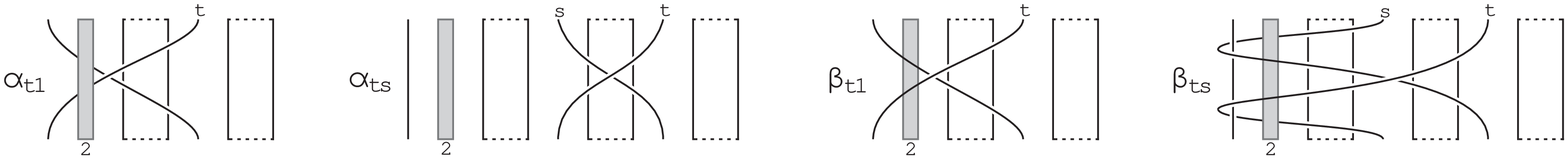}}
\end{center}
\vglue -0.3cm
\caption{Braid pictures for type~D new generators.}\label{fi:gend}
\centerline{\it Figure~\ref{fi:gend}: Diagrammes de tresse pour les nouveaux g\'en\'erateurs de type~D.}
\end{figure}

\ni We now consider type~D. The classical monoid~$\clas(D_n)$ for the Artin-Tits braid
group~$\gart(D_n)$ admits the presentation
\begin{equation}\label{eq:presclasd}
 \begin{aligned}
  \langle~\t_1^{},\s_1^{},\ldots,\s_{n-1}^{}:
   &~\s_1^{}\t_1^{}=\t_1^{}\s_1^{}~, ~\s_2^{}\t_1^{}\s_2^{}=\t_1^{}\s_2^{}\t_1^{}~,
   ~\s_i^{}\s_{i+1}^{}\s_i^{}=\s_{i+1}^{}\s_i^{}\s_{i+1}^{}~,~1\leq i\leq n-2\\
   &~\s_i^{}\s_j^{}=\s_j^{}\s_i^{}~,~1<i+1<j<n~,
   ~\t_1^{}\s_j^{}=\s_j^{}\t_1^{}~,~2<j<n~\rangle.
 \end{aligned}
\end{equation}
Let us introduce the following $n(n-1)$ new generators~:
\begin{equation}\label{eq:gend}
\begin{aligned}
\a_{ts}^{}&=(\s_{t-1}^{}\s_{t-2}^{}\cdots\s_{s+1}^{})
\s_s^{}(\s_{t-1}^{}\s_{t-2}^{}\cdots\s_{s+1}^{})^{-1}\qquad\hbox{for}\quad
n\geq t>s\geq 1,\\
\b_{t1}^{}&=(\s_{t-1}^{}\s_{t-2}^{}\cdots\s_2^{})
\t_1^{}(\s_{t-1}^{}\s_{t-2}^{}\cdots\s_2^{})^{-1}\qquad\hbox{for}\quad
n\geq t>1,\\
\b_{ts}^{}&=\a_{s1}^{-1}\b_{t1}^{}\a_{s1}^{}\qquad\hbox{for}\quad
n\geq t>s> 1.
\end{aligned}
\end{equation}
We use the pictures for the type D braids introduced by~Allcock in~\cite{all}, see~Figure~\ref{fi:gend}.

\begin{e-proposition}\label{presd} The dual braid monoid~$\dual(D_n)$ admits the presentation
\begin{equation}\label{eq:presd}
 \begin{aligned}
  \langle~\a_{ts},\b_{ts}:
   &~[\a_{ts},\a_{sr},\a_{tr}]~,~[\a_{ts},\b_{sr},\b_{tr}]~~\hbox{for}~~t>s>r,\\
   &~[\b_{ts},\a_{sr},\b_{tr}]~,~[\b_{tr},\a_{s1}]~,~[\b_{tr},\b_{s1}]~~\hbox{for}~~t>s>r>1,\\
   &~[\b_{ts},\b_{t1},\a_{s1}]~,~[\b_{ts},\a_{t1},\b_{s1}]~~\hbox{for}~~t>s>1,\\
   &~[\a_{ts},\a_{rq}]~,~[\a_{ts},\b_{rq}]~,~[\a_{tq},\a_{sr}]~,~[\b_{tq},\a_{sr}]~~\hbox{for}~~t>s>r>q,\\
   &~[\b_{ts},\a_{rq}]~~[\b_{tq},\b_{sr}]~~\hbox{for}~~t>s>r>q>1,\\
   &~[\a_{t1},\b_{t1}]~~\hbox{for}~~t>1~\rangle.
 \end{aligned}
\end{equation}
\end{e-proposition}

\begin{figure}[!hb]
\begin{center}
\scalebox{.90}{\includegraphics[width=17cm,height=4cm]{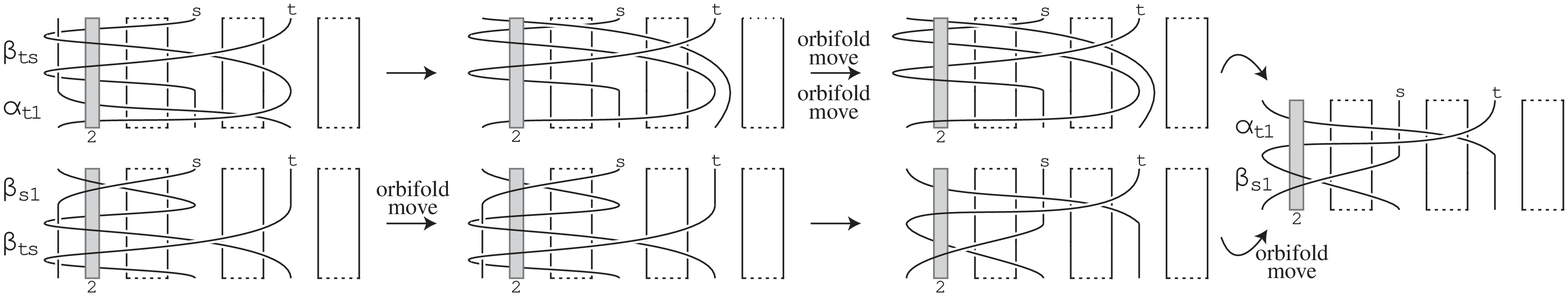}}
\end{center}
\vglue -0.3cm
\caption{Type~D braid isotopies.}\label{fi:disotopy}
\centerline{\it Figure~\ref{fi:disotopy}: Isotopies de tresses de type~D.}
\end{figure}

\begin{proof}(i) Applying the {\it
orbifold move} of~\cite{all}, the type~D braid isotopies  displayed in
Figure~\ref{fi:disotopy} give~$[\b_{ts},\a_{t1},\b_{s1}]$ for~$t>s>1$. The rest of this part of the proof is as
for~Proposition~\ref{presb}.

(ii) As for the type~B case, we have to complete Presentation~(\ref{eq:presd}), and the completion
algorithm gives~:
\vsm
$$\a_{tr}^{}\a_{rq}^{}\a_{ts}^{}=\a_{sq}^{}\a_{sr}^{}\a_{tq}^{},
\quad\a_{tr}^{}\a_{ts}^{}\b_{rq}^{}=\b_{sq}^{}\a_{sr}^{}\b_{tq}^{},
\quad\b_{ts}^{}\a_{tq}^{}\a_{sr}^{}=\b_{rq}^{}\b_{tr}^{}\b_{sq}^{},\vsm$$for~$t>s>r>q$,
\vsm
$$\a_{tq}^{}\a_{q1}^{}\a_{tr}^{}\a_{ts}^{}\b_{q1}^{}=\b_{sr}^{}\b_{s1}^{}\a_{rq}^{}\a_{s1}^{}\b_{tq}^{},
\quad\b_{tr}^{}\a_{rq}^{}\b_{ts}^{}=\a_{sq}^{}\a_{sr}^{}\b_{tq}^{},
\quad\b_{tr}^{}\a_{rq}^{}\a_{ts}^{}=\b_{sq}^{}\b_{sr}^{}\b_{tq}^{},\vsm$$for~$t>s>r>q>1$,\vsm\begin{equation*}
\begin{array}{rclrclrcl}
 \a_{tr}^{}\a_{r1}^{}\a_{ts}^{}\b_{r1}^{}\!\!&\!\!=\!\!&\!\!\b_{sr}^{}\b_{s1}^{}\a_{s1}^{}\b_{tr}^{},
 &\a_{t1}^{}\a_{ts}^{}\b_{r1}^{}\!\!&\!\!=\!\!&\!\!\b_{sr}^{}\a_{s1}^{}\b_{tr}^{},
 &\b_{ts}^{}\b_{t1}^{}\a_{sr}^{}\a_{t1}^{}\!\!&\!\!=\!\!&\!\!\a_{tr}^{}\a_{r1}^{}\a_{ts}^{}\b_{r1}^{},\\
 \b_{ts}^{}\b_{t1}^{}\a_{sr}^{}\!\!&\!\!=\!\!&\!\!\a_{r1}^{}\a_{s1}^{}\b_{tr}^{},
 &\b_{ts}^{}\b_{t1}^{}\a_{sr}^{}\a_{t1}^{}\!\!&\!\!=\!\!&\!\!\b_{sr}^{}\b_{s1}^{}\a_{s1}^{}\b_{tr}^{},
 &\b_{t1}^{}\a_{r1}^{}\a_{ts}^{}\!\!&\!\!=\!\!&\!\!\b_{sr}^{}\b_{s1}^{}\b_{tr}^{},\vsm\vsm 
\end{array}
\end{equation*}
for~$t>s>r>1$, and~$\a_{ts}^{}\a_{s1}^{}\b_{s1}^{}=\b_{ts}^{}\b_{t1}^{}\a_{t1}^{}$
for~$t>s>1$. Now, every relation in the complete presentation involves at most 4 distinct strands plus possibly
the first strand, so, as in the type~B case, checking local cube condition of the complete presentation given
for~$\dual(D_6)$---which is also immediate with a computer---is sufficient to prove global cube condition for~$\dual(D_n)$ for
every~$n$. Finally, $\d=\a_{n,n-1}\cdots\a_{2,1}\b_{2,1}$ is a Garside element, whose image in the Coxeter
group is a Coxeter element.
\end{proof}

\section{Combinatorics of the dual Garside structures}\label{comb}

\ni Birman, Ko \&~Lee showed in~\cite{bkl} that the simple elements of the dual braid monoid~$\dual(A_{n-1})$ are in one-to-one correspondence
with the non-crossing partitions of the integer~$n$ (see also~\cite{bdm}). An analogous approach allows us to prove that the simple
elements of~$\dual(B_n)$ and~$\dual(D_n)$ are in bijection with the corresponding Reiner's non-crossing partitions of~\cite{rei}.
Table~\ref{ta:table} gives the number of simple elements for the dual braid monoids, gathering theoritical results for~A, B, D
and~I$_2$ types and computational results---using the package CHEVIE of~GAP~\cite{gap}---for exceptional types (the computation
for~E$_8$ first appeared in~\cite{bes}).

\begin{center}
\begin{table}[!ht]
\begin{tabular}{|c|c|c|c|c|c|c|c|c|c|c|}
\hline
type&A$_n$&B$_n$&D$_n$&$\!$H$_3\!$&F$_4$&H$_4$&E$_6$&E$_7$&E$_8$&$\!\!$I$_2(m)\!\!$\\
\hline
&&&&&&&&&&\\
\!classical\!
&$(n+1)!$&$2^nn!$&$2^{n-1}n!$
&$\!\!\!\!120\!\!\!\!$&$\!\!1152\!\!$&$\!\!14400\!\!$
&$\!\!51840\!\!$&$\!\!2903040\!\!$&$\!\!696729600\!\!$
&$2m$\\
&&&&&&&&&&\\
\hline
&&&&&&&&&&\\
dual
&$\!{1\over n+2}{2n+2\choose n+1}\!$&${2n\choose n}$&$\!{2n\choose n}\!\!-\!\!{2n-2\choose n-1}\!$
&$\!32\!$&$\!105\!$&$\!280\!$
&$\!833\!$&$\!4160\!$&$\!25080\!$
&$\!\!m+2\!\!$\\
&&&&&&&&&&\\
\hline
\end{tabular}
\caption{The number of simple elements in classical and dual braid monoids.}\label{ta:table}
\centerline{\it Table~\ref{ta:table}: Le nombre d'\'el\'ements simples dans les mono\"\i des de tresses classiques
et duaux.}
\end{table}
\end{center}


%
\end{document}